\documentclass[11pt]{paper}

\usepackage{color}
\usepackage{mathptmx}       
\usepackage{helvet}         
\usepackage{courier}        
\usepackage{type1cm}        
\usepackage{makeidx}         
\usepackage{graphicx}        
\usepackage{multicol}        
\usepackage[bottom]{footmisc}
\usepackage{amsmath}
\usepackage{amsfonts}
\usepackage{amssymb}
\usepackage{amsthm}
\newtheorem{theorem}{Theorem}

\newtheorem{lemma}[theorem]{Lemma}

\newtheorem{proposition}[theorem]{Proposition}

\begin{document}

\title{From CLE($\kappa$) to SLE($\kappa,\rho$)'s}
\author{Wendelin Werner \and Hao Wu\thanks{H.W. work is funded by the Fondation CFM JP Aguilar pour la recherche. The authors acknowledge also the support and/or hospitality of the ANR grant BLAN-MAC2, of the Einstein Foundation, of TU Berlin and FIM at ETH Z\"urich}
}
\date {Universit\'e Paris-Sud}

%
%
\maketitle

\abstract{We show how to connect together the loops of a simple Conformal Loop Ensemble (CLE) in order to construct samples of
chordal SLE$_\kappa$ processes and their SLE$_{\kappa}(\rho)$ variants, and we discuss some consequences of this construction.}

\newcommand{\eps}{\epsilon}
\newcommand{\ov}{\overline}
\newcommand{\U}{\mathbb{U}}
\newcommand{\T}{\mathbb{T}}
\newcommand{\HH}{\mathbb{H}}
\newcommand{\LA}{\mathcal{A}}
\newcommand{\LF}{\mathcal{F}}
\newcommand{\LK}{\mathcal{K}}
\newcommand{\LE}{\mathcal{E}}
\newcommand{\LL}{\mathcal{L}}
\newcommand{\R}{\mathbb{R}}
\newcommand{\C}{\mathbb{C}}
\newcommand{\N}{\mathbb{N}}
\newcommand{\Z}{\mathbb{Z}}
\newcommand{\E}{\mathbb{E}}
\newcommand{\PP}{\mathbb{P}}
\newcommand{\QQ}{\mathbb{Q}}
\newcommand{\MR}{MR}

\noindent
\section{Introduction}

The goal of the present paper is to derive ways to construct samples (chordal) SLE curves (or the related SLE$_{\kappa}(\rho)$ curves) out of the sample of a Conformal Loop Ensemble (CLE), using additional Brownian paths (or so-called restriction measure samples). In order to properly state a first version of our result, we need to briefly informally recall the definition of these three objects: SLE, CLE and the restriction measures.

\begin{itemize}
\item
Recall that a chordal SLE (for Schramm-Loewner Evolution) in a simply
connected domain $D$ is a random curve that is joining two prescribed boundary points $a$ and $b$ of $D$. These curves have been first defined by Oded Schramm in 1999 \cite {OS2000}, who conjectured (and this conjecture was since then proved in several important cases) that they should be the scaling limit of particular random curves in two-dimensional critical statistical physics models when the mesh of the lattice goes to $0$. More precisely, one has typically to consider the statistical physics model in a discrete lattice-approximation of $D$, with well-chosen boundary conditions, where (lattice-approximations of) the points $a$ and $b$ play a special role. When $\kappa \le 4$, these SLE$_\kappa$ curves are random
simple continuous curves that join $a$ to $b$ with fractal dimension is $1+ \kappa/8$ (see for instance \cite {GL2005} and the references therein).

\item
CLEs (for Conformal Loop Ensembles) are closely related objects. A CLE is a random family of loops that is defined in a simply connected domain $D$.
In the present paper, we will only discuss the CLEs  that consist of simple loops. There are various equivalent definitions and constructions of these simple CLEs -- see for instance the discussion in \cite {SW2010}. More precisely, one CLE sample is a collection of countably many disjoint simple loops in $D$, and it is conjectured to correspond to the scaling limit of the collection of all discrete (but macroscopic) interfaces in the corresponding lattice model from statistical physics. Here, the boundary conditions are ``uniform'' and involve no special marked points on the boundary of $D$ (as opposed to the definition of chordal SLE that requires to choose the boundary points $a$ and $b$).
It is proved in \cite {SW2010} that there is exactly a one-dimensional family of simple CLEs, that is indexed by $\kappa \in (8/3, 4]$. Then, in a CLE$_\kappa$ sample, the loops all locally look like SLE$_\kappa$ type curves (and have fractal dimension $1+ \kappa/8$). Note also that, even if any two loops are disjoint in CLE$_\kappa$ sample, the Lebesgue measure of the set of points that are surrounded by no loop is almost surely $0$. This is therefore a random Cantor-like set, sometimes called the CLE carpet (its fractal dimension is actually proved in \cite {SSW,NW2011} to be equal to $1 + (2 / \kappa) + 3 \kappa/32 \in [15/8, 2)$). In the present paper, we will only discuss the CLEs for $\kappa \le 4$, that consist of simple disjoint loops (there exists other CLEs for $\kappa \in (4, 8]$).
\item
When $a$ and $b$ are two boundary points of a simply connected domain $D$ as before, it is possible to define random simple curves from $a$ to $b$ that posess a certain ``one-sided restriction'' property, that is defined and discussed in \cite {LSW2003}. There is in fact a one-dimensional family of such random curves, that is parametrized by its restriction exponent, which can take any positive real value $\alpha$. All these random restriction curves can be viewed as boundaries of certain Brownian-type paths (or like SLE$_{8/3}$ curves). In particular, they all almost surely have a Hausdorff dimension that is equal to $4/3$.

\end{itemize}

\begin{figure}[ht!]
\begin{center}
\includegraphics[width=\textwidth]{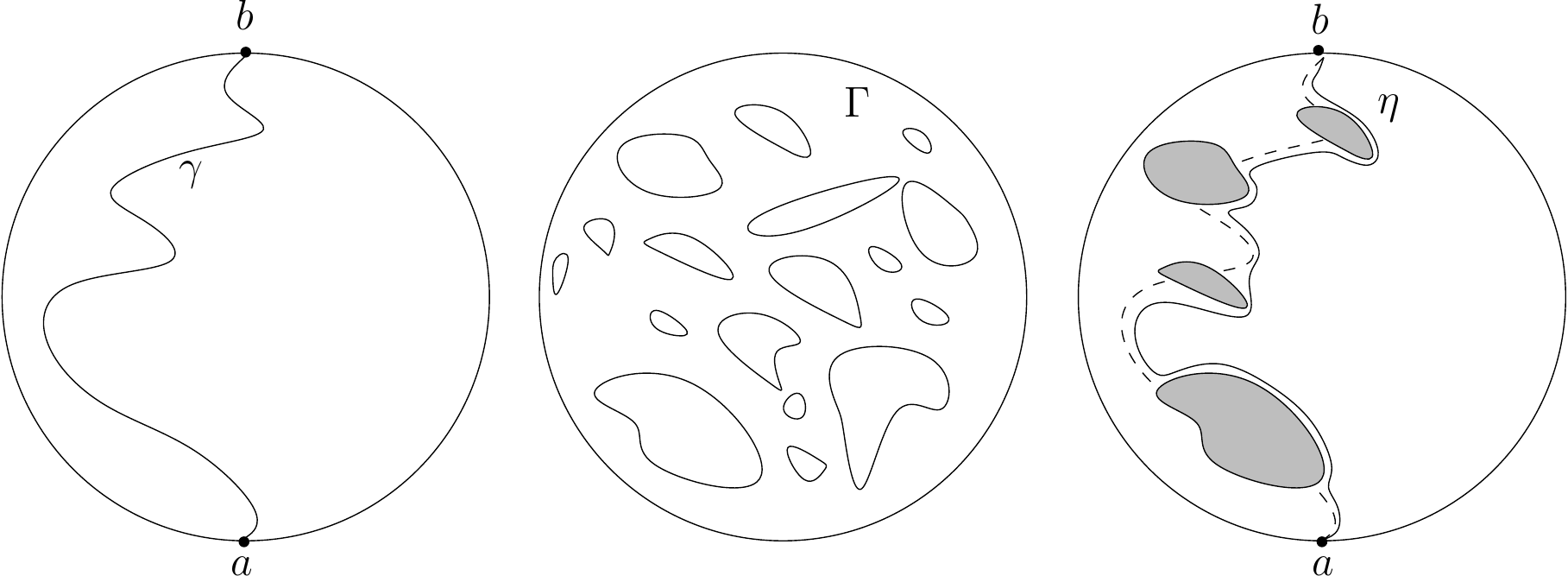}
\end{center}
\caption{\label{fig::construction_eta}Construction of $\eta$ out of $\gamma$ and $\Gamma.$}
\end{figure}

Let us now state the main result that we prove in the present paper: Define independently, in a simply connected domain $D$ with two marked boundary points $a$ and $b$, the following two random objects: A CLE$_\kappa$ (for some $\kappa \in (8/3, 4]$) that we call $\Gamma$ and a one-sided restriction path $\gamma$ from $a$ to $b$, with restriction exponent $\alpha$.
Finally, we define the set obtained by attaching to $\gamma$ all the loops of $\Gamma$ that it intersects. Then, we define the right-most boundary of this set. This turns out to be again a simple curve from $a$ to $b$ in $D$ that we call $\eta$ (see Figure \ref{fig::construction_eta}). Note that in order to construct $\eta$, it is enough to know $\gamma$ and the outermost loops of $\Gamma$.

\begin {theorem}
\label {main}
When $\kappa \in (8/3, 4]$ and $\alpha = (6-\kappa)/(2\kappa)$, then $\eta$ is a chordal SLE$_\kappa$ from $a$ to $b$ in $D$.
\end {theorem}

In fact, for a given $\kappa$, the other choices of $\alpha >0$ give rise to variants of SLE$_\kappa$, the so-called SLE$_{\kappa}(\rho)$ curves,  where $\rho$ is related to $\kappa$ and $\alpha$ by the relation  $\alpha= ( \rho+2)(\rho+ 6 - \kappa)/ (4 \kappa)$. We will state this generalization of Theorem \ref{main} in the next section, after having properly introduced these SLE$_{\kappa}(\rho)$ processes.

\medbreak

To illustrate Theorem \ref {main}, let us give the following example for $\kappa=3$, which corresponds to the scaling limit of the critical Ising model (see \cite {CS2012,CDHKS}). Consider a CLE$_3$ $\Gamma$ in $D$ which is the (soon-to-be proved) scaling limit of the collection of outermost critical Ising model  ``$-$ cluster'' boundaries, when one considers the model with uniformly ``$+$ boundary conditions''. On the other hand, consider now the scaling limit of the critical Ising model with mixed boundary conditions, $+$ between $a$ and $b$ (anti-clockwise) and $-$ between $b$ and $a$. This model defines loops as before, as well as the additional $\pm$ interface $\eta$ joining $a$ and $b$, which turns out to be a SLE$_3$ path (see \cite {CDHKS}).
Now, our result shows that in order to construct a sample of $\eta$,
one possibility is to take the right boundary of the union of a restriction measure with exponent $1/2$ together with all the loops in $\Gamma$ that it intersects. It gives a way to see the ``effect'' of changing the boundary conditions (note that there are natural ways to couple the discrete Ising model with mixed boundary conditions to the model with uniform boundary conditions, it would be interesting to compare them with this
coupling in the scaling limit).

\medbreak

We would like to make a few comments:
\begin {enumerate}
 \item
 It is proved in \cite {SW2010} that CLEs can be constructed as outer boundaries of clusters of Poissonian clouds of Brownian loops in $D$ (the ``Brownian loop-soups'' introduced in \cite {LW2004}) with intensity $c(\kappa)$. Hence, together with the construction of the restriction measure via clouds of Brownian excursions or reflected Brownian motions, this provides a ``completely Brownian'' construction of all these chordal SLE$_\kappa$ curves and their SLE$_\kappa(\rho)$ variants. This result was in fact announced in \cite {WW2003}, so that -- combined with \cite {SW2010} -- the present paper eventually completes the proof of that (not so recent) research announcement.

 \item This Brownian construction of SLE$_\kappa(\rho)$ paths turn out to be particularly useful and handy, when one has to derive ``second moment estimates'' for these SLE curves. We will illustrate this in the final section of the present paper by giving a short self-contained derivation of the Hausdorff dimension of the intersection of SLE$_\kappa (\rho)$ (in the upper half-plane) with the real line.

 \item
 A direct by-product of this construction of these chordal SLE$_\kappa$ curves and their variants is that they are ``reversible'' simple paths (for instance, the SLE from $a$ to $b$ in $D$ is a simple path has the same law as the SLE from $b$ to $a$ modulo reparametrization -- in the case of SLE$_\kappa (\rho)$ the statement is also clear, but the reversed SLE$_\kappa (\rho)$ is then pushed/attracted from its right). This provides an alternative proof to the reversibility of these SLE$_{\kappa}(\rho)$ curves that has been obtained thanks to their relation with the Gaussian Free Field in \cite {MS2012b} (see also \cite {zhan0,zhan,JD2009} for earlier proofs of this result in the case $\rho=0$ and then when the SLE$_\kappa (\rho)$ curves do not hit the boundary of the domain i.e. when $\rho \ge (\kappa -4)/2$). Note however that our approach does not yield any result for $\kappa \notin [ 8/3, 4]$.

 \item
 The construction of the restriction measure via Poisson point processes of Brownian excursions, as explained in \cite {WW2005}, together with that of the CLE's via loop-soups, make it possible to define simultaneously in a fairly natural and ``ordered way'' (see the comments after the statement of Theorem \ref {main2}), on a single probability space, all these SLE$_{\kappa}(\rho)$'s in $D$ from $a$ to $b$, for all boundary points $a$ and $b$, and for all $\kappa \in (8/3, 4]$ and all $\rho > -2$.
 This is of course reminiscent of the definitions of SLE$_{\kappa}(\rho)$ processes within a Gaussian Free Field \cite {MS2012a}. It is interesting to see the similarities and differences between these two constructions.

 \end {enumerate}

\section{Preliminaries}
In this section, we will recall in a little more detail some definitions, notations and facts, and point to appropriate references for background. We then state our main result, Theorem \ref {main2} and make a couple of remarks.

\subsection{Conformal restriction property}
We first recall the definition and the basic properties of the paths satisfying conformal restriction (almost all the results that we shall describe have been derived in \cite {LSW2003}, a survey as well as the construction
of restriction samples from Brownian excursions can be found in \cite {WW2005}).

Here and throughout the paper, we denote the upper half of the complex plane $\C$ by
$\HH:=\{x+iy: x\in\R, y>0\}$.
Let $\LA$ be the set of all bounded closed $A\subset \overline{\HH}$ such that $\R_- \cap  A = \emptyset$ and $H_A := \HH\setminus A$ is simply connected.

For $A\in\LA,$ we define $\Phi_A$ to be the unique conformal map from $H_A$ onto $\HH$ such that $\Phi(0)= 0$ and $\Phi_A(z) \sim z $ as $z\to\infty$.

\begin{figure}[ht!]
\begin{center}
\includegraphics[width=0.7\textwidth]{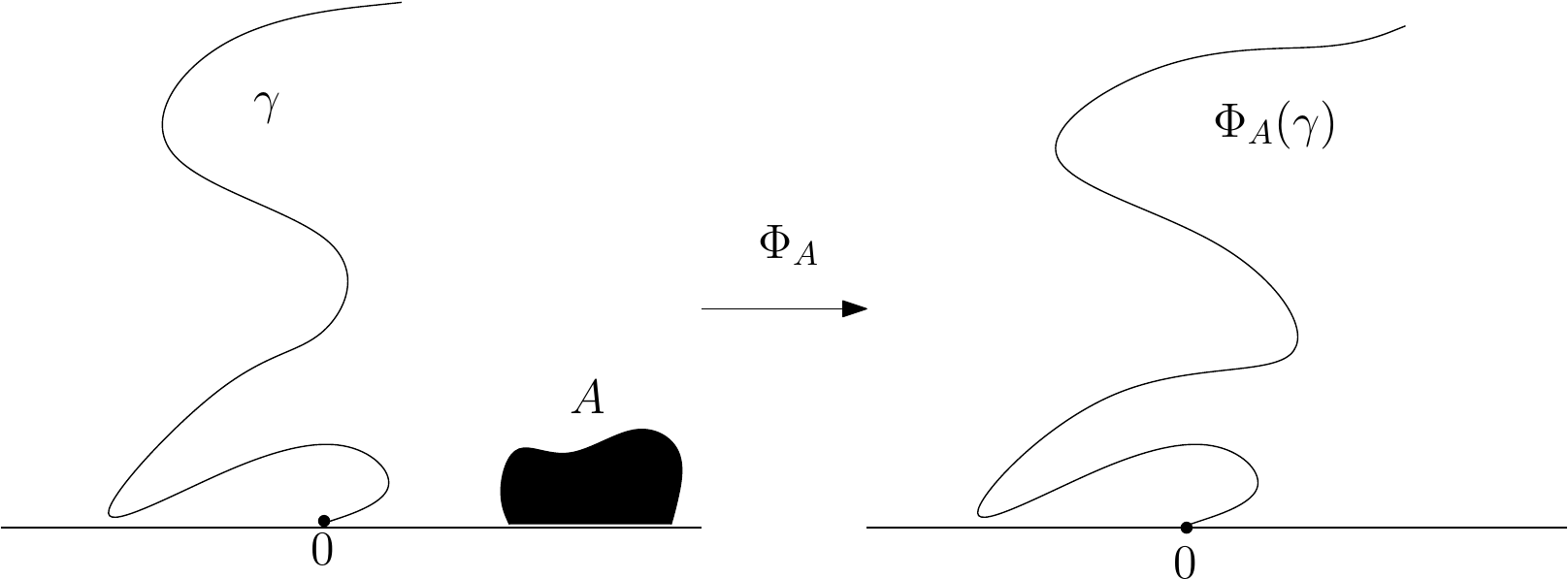}
\end{center}
\caption{\label{fig::conformal_restriction} The law of $\Phi_A(\gamma)$ conditionally on $\gamma\cap A=\emptyset$ has the same law as $\gamma$ itself.}
\end{figure}

We say that a random curve $\gamma$ from $0$ to infinity in $\overline \HH$ does satisfy one-sided conformal restriction (to the right), if for any $A$, the law of $\Phi_A (\gamma)$ conditionally on $\gamma \cap A = \emptyset$ is in fact identical to the law of $\gamma$ itself (see Figure \ref{fig::conformal_restriction}).

It turns out that if this is the case, then there exists some non-negative $\alpha$ such that for all $A \in \LA$,
\begin {equation}
 \label {alpha}
P ( \gamma \cap A=\emptyset)=\Phi_A'(0)^{\alpha}.
\end {equation}
Conversely, for all non-negative $\alpha$, there exists exactly one distribution for $\gamma$ that fulfils (\ref{alpha}) for all $A \in \LA$. We call $\gamma$ an one-sided restriction sample of exponent $\alpha.$
There exist several equivalent constructions of $\gamma$:
\begin {itemize}
 \item As the right boundary of a certain Brownian motion from $0$ to $\infty$, reflected on $(-\infty, 0]$ with a certain reflection angle $\theta (\alpha)$ and conditioned not to intersect $[0, \infty)$, see \cite {LSW2003}.
 \item As the right boundary of Poissonian cloud of Brownian excursions from $(-\infty, 0]$ in $\HH$ (so it is the right boundary of the countable union of Brownian paths that start and end on the negative half-line, see \cite {WW2005}).
 \item As an SLE$_{8/3}(\rho)$ curve for some $\rho > -2$ (these processes will be defined in the next subsection), see \cite {LSW2003} for the relation between $\alpha$ and $\rho$. Note that this approach enables to show that $\gamma$ does hit the negative half-line if and only $\alpha < 1/3$.
\end {itemize}

We can note that the limiting case $\alpha = 0$ corresponds to the case where $\gamma$ is the negative half-line, whereas the case $\alpha = 5/8$ corresponds to $\rho = 0$ i.e. to the SLE$_{8/3}$ curve itself, which is left-right symmetric. Furthermore, the second construction shows immediately that for $\alpha < \alpha'$, it is possible to couple the corresponding restriction curves in such a way that $\gamma'$ stays ``to the right'' of $\gamma$ (with obvious notation). In other words, the larger $\alpha$ is, the more the restriction sample is ``repelled'' from the negative half-line.

In fact, we will be only using the second description in the present paper (and we will actually recall in Subsection \ref {Smain} why this indeed constructs a random simple curve $\gamma$).

\subsection{SLE$_{\kappa}(\rho)$ process}
\label {2.2}

The SLE$_{\kappa}(\rho)$ processes are natural variants of SLE$_\kappa$ processes that have been first introduced in \cite {LSW2003}. Recall first  that
the SLE$_\kappa$ curves for $\kappa \le 4$ are random simple continuous curves $\eta$ from $0$ to $\infty$ in $\HH$ that possess the following properties:
\begin {itemize}
 \item The law of $\eta$ is scale-invariant: For any positive $\lambda$, the traces of $\eta$ and of $\lambda \eta$ have the same law.
 \item Let us suppose that $\eta$ is parametrized by its half-plane capacity (i.e., for any $t$, the conformal map $g_t$ from $\HH \setminus \eta [0,t]$ onto
 $\HH$ such that $g_t (z) \sim z + o(1)$ when $z \to \infty$ in fact satisfies $g_t(z) - z \sim 2t/z $). For any positive time $t$, the distribution of $g_t ( \eta[t, \infty)) - g_t (\eta_t)$ is identical to the distribution of $\eta$ itself.
\end {itemize}
In fact, the SLE$_\kappa$ curves are the only random curves with this property, which is what led Oded Schramm to the definition of these curves, that involves the Loewner differential equation describing growing hulls, where one chooses the driving function to be a one-dimensional Brownian motion (see \cite {OS2000}).

There exist variants of the SLE$_\kappa$ curves that involve additional marked boundary points, and that are called the 
SLE$_{\kappa}(\rho_1, \ldots, \rho_L)$ processes. Let us now describe the SLE$_{\kappa}(\rho)$ processes that involve exactly one additional marked boundary point (see \cite {LSW2003,JD2005}).

\begin{figure}[ht!]
\begin{center}
\includegraphics[width=0.9\textwidth]{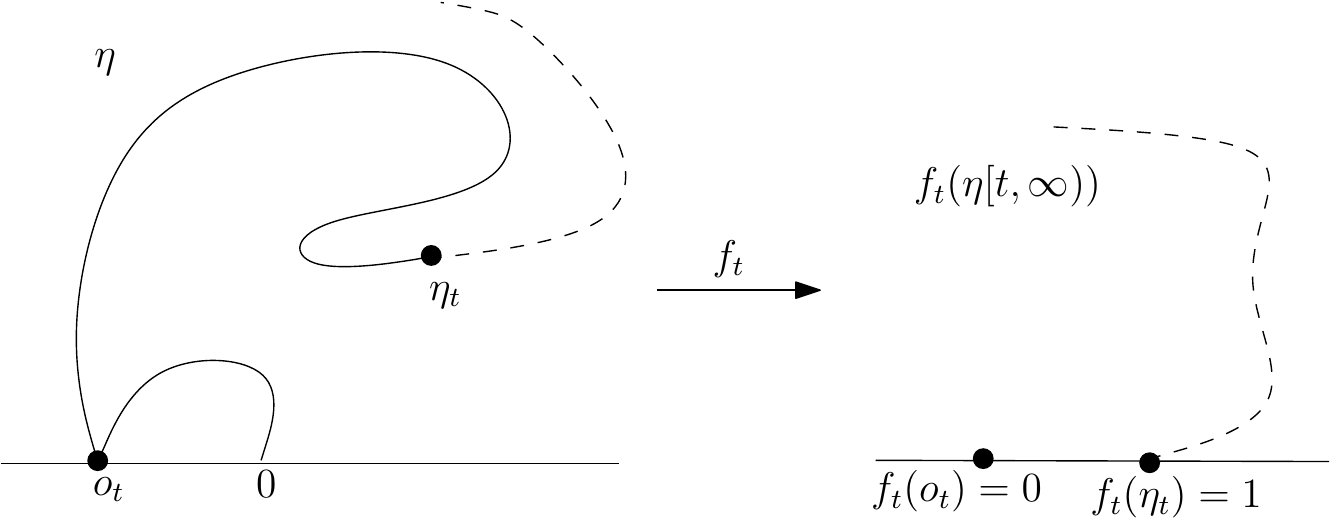}
\end{center}
\caption{\label{fig::SLE_kappa_rho} $f_t(\eta[t,\infty))$ is independent of $\eta[0,t].$}
\end{figure}

It turns out that they can also be characterized by a couple of properties. Let us now state the characterization that will be handy for our purposes:
Suppose that the following four properties hold:
\begin {itemize}
\item $\eta$ is a random simple curve from $0$ to $\infty$ in $\overline \HH$.
\item The law of $\eta$ is scale-invariant: For any positive $\lambda$, the traces of $\lambda \eta$ and $\eta$ are identically distributed.
\item $\eta \cap (0, \infty) = \emptyset$ and the Lebesgue measure of $\eta \cap (-\infty, 0]$ is almost surely equal to $0$. Mind however that it is possible (and it will happen in a number of cases) that $\eta$ hits the negative half-line.
\item Suppose that $\eta$ is parametrized by half-plane capacity as before. For any positive time $t$, define $H_t$ as the unbounded connected component of $\HH\setminus\eta[0,t]$ (if $\eta$ intersects the negative half-line, it happens that $H_t \not= \HH \setminus \eta [0,t]$) and $o_t$ as the left-most point of the intersection $\eta[0,t]\cap \R_-.$ Let $f_t$ be the unique conformal map from $H_t$ onto $\HH$ such that sends the triplet $(o_t,\eta_t,\infty)$ onto $(0,1,\infty).$ Then, the distribution of $f_t (\eta[t, \infty))$ is independent of $t$ (and of $\eta [0,t]$) (see Figure \ref{fig::SLE_kappa_rho}).
\end {itemize}

Then, $\eta$ is necessarily a SLE$_{\kappa}(\rho)$ for some $\kappa \in (0, 4]$ and $\rho > -2$ (mind that the fact that this SLE$_\kappa (\rho)$ is almost surely a simple curve is then part of the conclusion; in fact in the present paper, we will never use the {\sl a priori} fact that the SLE$_\kappa (\rho)$ processes are continuous simple paths).

This is very easy to see, using the Loewner chain description of the random simple curve $\eta$.
If one parametrizes the curve $\eta$ by its half-plane capacity (which is possible because the its capacity is increasing continuously -- this is due to the third property) and defines the usual conformal map $g_t$ from $H_t$ onto $\HH$ normalized by $g_t (z) = z + o(1)$ near infinity, then one can define
$$ W_t = g_t (\eta_t), O_t = g_t (o_t).$$
One observes that $X_t := W_t - O_t$ is a Markov process with the Brownian scaling property i.e., a multiple of a Bessel process.
More precisely, one can first note that the first two items imply that for any given $t_0> 0$, $\eta_{t_0} \notin (- \infty, 0)$ and therefore $u:=X_{t_0} \not= 0$. The final property then implies readily that the law of
$((X_{t_0 + tu^2} -u) / u, t \ge 0)$ is independent of $(X_t, t \le t_0)$. From this, it follows that at least up to the first time after $t_0$ at which $X$ hits the origin, it does behave like a Bessel process. Then, one can notice that $X$ is instantaneously reflecting away from $0$ because the Lebesgue measure of the set of times at which it is at the origin is almost surely equal to $0$. Hence, one gets that $X$ is the multiple of some reflected Bessel process of positive dimension (see \cite {RY1999} for background on Bessel processes). From this, one can then recover the process $t \mapsto O_t$ (because of the Loewner equation $dO_t = 2dt / (O_t - W_t)$ when $X_t \not= 0$) and finally $t \mapsto W_t$. In particular, we get that
$$ dW_t = \sqrt {\kappa} dB_t + \frac {\rho}{W_t -O_t} dt $$
for some $\rho > -2$ and $\kappa \le 4$ (the fact that $\rho > -2$ is a consequence of the fact that the dimension of the Bessel process $X$ is positive; $\kappa \le 4$ is due to the fact that $\eta$ does not hit the positive half-line). This characterizes the law of $\eta$, which is then called the SLE$_{\kappa}(\rho)$.

Actually, it is possible to remove some items from this characterization of SLE$_{\kappa}(\rho)$ curves; the first three items are slightly redundant, but since we do get these properties for free in our setting, the present presentation will be sufficient for our purposes (see for instance \cite {SchrammW2005, MS2012b} for a more general characterization).

Note that the SLE$_{\kappa}(\rho)$ processes touch the negative half-line if and only if $\rho < (\kappa/2) - 2$ (as this corresponds to the fact that the Bessel process $(W_t - O_t )/ \sqrt {\kappa}$ has dimension smaller than $2$).

Let us point out that it is possible to make sense also of SLE$_{\kappa}(\rho)$ processes for some values of $\rho \le -2$ by introducing either a symmetrization or a compensation procedure (see \cite {JD2005,Sheffield,WW}), some of which are very closely related to CLEs as well, but we will not discuss such generalized SLE$_{\kappa} (\rho)$'s  in the present paper.

\subsection {Simple CLEs}

Let us now briefly recall some features of the Conformal Loop Ensembles for $\kappa \in (8/3, 4]$ -- we refer to \cite {SW2010} for details (and the proofs) of these statements.
A CLE is a collection $\Gamma$ of non-nested disjoint simple loops $(\gamma_j, j \in  J)$ in $\HH$ that possesses a particular conformal restriction property. In fact, this property that we will now recall, does characterize these simple CLEs:
\begin {itemize}
\item For any M\"obius transformation $\Phi$ of $\HH$ onto itself, the laws of $\Gamma$ and  $\Phi(\Gamma)$ are the same. This makes it possible to define, for any simply connected domain $D$ (that is not the entire plane -- and can therefore be viewed as the conformal image of $\HH$ via some map $\tilde \Phi$), the law of the CLE in $D$ as the distribution of $\tilde \Phi (\Gamma)$ (because this distribution does then not depend on the actual choice of conformal map $\tilde \Phi$ from $\HH$ onto $D$).
\item For any simply connected domain $H \subset \HH$, define the set $\tilde H = \tilde H ( H, \Gamma)$ obtained by removing from $H$ all the loops (and their interiors) of $\Gamma$ that do not entirely lie in $H$. Then, conditionally on $\tilde H$, and for each connected component $U$ of $\tilde H$, the law of those loops of $\Gamma$ that do stay in $U$ is exactly that of a CLE in $U$.
\end {itemize}

It turns out that the loops in a given CLE are SLE$_\kappa$ type loops for some value of $\kappa \in (8/3, 4]$ (and they look locally like SLE$_\kappa$ curves).  In fact for each such value of $\kappa$, there exists exactly one CLE distribution that has SLE$_\kappa$ type loops.
As explained in \cite {SW2010}, a construction of these particular families of loops can be given in terms of outermost boundaries of clusters of the Brownian loops in a Brownian loop-soup with subcritical intensity $c$ (and each value of $c$ corresponds to a value of $\kappa$).

\subsection {Main Statement}
\label {Smain}
We can now state our main Theorem, that generalizes Theorem \ref {main}: Suppose that $\kappa \in (8/3, 4]$ is fixed (and it will remained fixed throughout the rest of the paper) and consider a CLE$_\kappa$ in the upper half-plane.
Independently, sample a restriction curve $\gamma$ from $0$ to $\infty$ in $\HH$ with positive exponent $\alpha$, and define $\eta$ out of the CLE and $\gamma$ just as in Theorem \ref {main}. Let $\tilde \rho := \tilde \rho ( \kappa, \alpha)$ denote the unique real in $(-2, \infty)$ such that
$$\alpha = \frac {(\tilde \rho + 2) (\tilde \rho+ 6 - \kappa)}{4\kappa}$$ (we will use this notation throughout the paper).
 Then:

\begin {theorem}
 \label {main2}
 The curve $\eta$ is a random simple curve which is an SLE$_{\kappa}(\tilde\rho)$.
\end {theorem}

Note that for a fixed $\kappa \in (8/3, 4]$, the function $\alpha \mapsto \tilde \rho$ is indeed an increasing bijection from $(0, \infty)$ onto $(-2, \infty)$.
The limiting case $\rho = -2$ in fact can be interpreted as corresponding to the case where both $\gamma$ and $\eta$ are the negative half-line.
Similarly, in the limiting case $\kappa=8/3$, where the CLE is in fact empty, then Theorem \ref {main2} corresponds to the description of $\gamma$ itself as an SLE$_{8/3}(\rho)$ curve.

\medbreak

Note that this construction shows that it is possible to couple an SLE$_\kappa (\rho)$ with an SLE$_{\kappa'}(\rho')$ in such a way that the former is almost surely  ``to the left'' of the latter, when $8/3 < \kappa \le \kappa' \le 4 $ and $\rho$ and $\rho'$ are chosen in such a way that
$$(\rho + 2) (\frac {\rho+ 6} \kappa  - 1 ) \le (\rho' + 2) (\frac {\rho'+ 6} {\kappa'}  - 1 ).$$
For example, an SLE$_\kappa (\rho)$ can be chosen to be to the left of an SLE$_\kappa (\rho')$ for $\rho \le \rho'$. Or an SLE$_3$ can be coupled to an SLE$_4(2 \sqrt{2}-2)$ in such a way that it remains almost surely to its left. Such facts are seemingly difficult to derive directly from the Loewner equation definitions of these paths.

Similarly, it also shows that it is possible to couple an SLE$_\kappa (\rho)$ from $0$ to $\infty$ with another SLE$_\kappa (\rho)$ from $1$ to $\infty$, in such a way that the latter stays to the ``right'' of the former.

\medbreak

Let us recall that the definition of SLE$_{\kappa}(\rho)$ processes can be generalized to more than one marked boundary point. For instance, if one considers
$x_1< \ldots < x_n \le 0 \le x_1' < x_2' < \ldots < x_l'$, it is possible to define a SLE$_{\kappa}(\rho_1, \ldots, \rho_n ; \rho_1', \ldots , \rho_l')$ from $0$ to infinity in $\HH$, with marked boundary points $x_1, \ldots, x_l'$ with corresponding weights. Several of these processes have also an interpretation in terms of conditioned SLE$_{\kappa}(\rho)$ processes (where the conditioning involves non-intersection with additional restriction samples) -- see \cite {WW2004}, so that they can also be interpreted via a CLE and restriction measures.

\medbreak

Let us now immediately explain why $\eta$ is necessarily almost surely a continuous curve from $0$ to $\infty$ in $\overline \HH$.
Let us first map all items (the CLE loops and the restriction sample) onto the unit disc, via the Moebius map $\Phi$ that maps $0$, $i$ and $\infty$ respectively onto
$-1$, $0$ and $1$, and write $\tilde \Gamma =\Phi (\Gamma)$, $\tilde \eta= \Phi (\eta)$ and $\tilde \gamma = \Phi (\gamma)$.

Let us first recall from \cite {SW2010} that $\tilde \Gamma$ consists of a countable family of disjoint simple loops such that for any $\eps>0$, there exist only finitely many loops of diameter greater than $\eps$.
Let us note that $\tilde \gamma$ is almost surely a continuous curve from $-1$ to $1$ in the closed unit disc. One simple way to check this (but other justifications are possible) is to use the construction of $\tilde \gamma$ as the bottom boundary of the union of countably many excursions away from the top half-circle. More precisely, for each excursion $e$ in this Poisson point process, one can define the loop $l(e)$ obtained by adding to this excursion the arc of the top half-circle that joins the endpoints of $e$. Then, one can construct a continuous path $\lambda$ from $-1$ to $1$ by moving from $-1$ to $1$ on this top arc, and attaching all these loops $l(e)$ in the order in which one meets them. As almost surely, for any $\eps>0$, there are only finitely many loops $l(e)$ of diameter greater than $\eps$, there is a way to parametrize $\lambda$ as a continuous function from $[0,1]$ into the closed disk. We then complete $\lambda$ into a loop by adding the bottom half-circle. Then, we can interpret $\tilde \gamma$ as part of the boundary of a connected component of the
complement of a continuous loop in the plane: It is therefore necessarily a continuous curve and it is easy to check that it is self-avoiding (because the Brownian excursions have no double cut-points).

We have detailed the previous argument, because it can be repeated in almost identical terms to explain why $\tilde \eta$ is a simple curve: We now move along $\tilde \gamma$ and attach the loops of $\tilde \Gamma$ that it encounters, in their order of appearance. By an appropriate time-change, we can ensure that the obtained path that joins $-1$ to $1$ in the closed disk is a continuous curve from $[0,1]$ into the closed unit disk. Then, just as above, we complete this curve into a loop by adding the bottom half-circle, and note that $\tilde \eta$ is a continuous curve from $-1$ to $1$. It is then easy to conclude that it is self-avoiding, because almost surely, $\tilde \gamma$ does never hit a loop
of $\tilde \Gamma$ at just one single point (this is due to the Markov property of Brownian motion: If one samples first the CLE and then the Brownian excursions that are used to construct $\gamma$, almost surely, a Brownian excursion will actually enter the inside of each individual loop of $\Gamma$ that it hits).

\section {Identification of $\rho$}

The proof of Theorem \ref {main2} consists of the following two steps.
\begin {lemma}
 \label {firstpart}
 The random simple curve $\eta$ is an SLE$_{\kappa}(\rho)$ curve for some  $\rho > -2$.
\end {lemma}
\begin {lemma}
 \label {secondpart}
 If $\eta$ is an SLE$_{\kappa}(\rho)$  for some $\rho > -2$, then necessarily $\rho = \tilde \rho (\kappa, \alpha)$.
\end {lemma}

The proof of Lemma \ref {firstpart} will be achieved in the next section by proving that it satisfies all the properties that characterize these curves (and that we have recalled in the previous subsection), which is the most demanding part of the paper.
In the present section, we will prove Lemma \ref {secondpart}. These ideas were already very briefly sketched in \cite {WW2003}.

\medbreak

Let us build on the loop-soup cluster construction of the CLE$_\kappa$ as established in \cite {SW2010}. We therefore consider a Poisson point process of Brownian loops (as defined in \cite {LW2004}) in the upper-half plane with intensity $c(\kappa) \in (0,1]$
with
$$c(\kappa)=\frac{(3\kappa-8)(6-\kappa)}{2\kappa}.$$
Then, we construct the CLE$_\kappa$ as the collection of all outermost boundaries of clusters of Brownian loops (here, we say that two loops $l,l'$ in the loop-soup are in the same cluster of loops if one find a finite chain of loops $l_0,...,l_n$ in the loop-soup such that $l_0=l, l_n=l,$ and $l_j\cap l_{j-1}\neq\emptyset$ for $j\in\{1,...,n\}$), as explained in \cite {SW2010}.

We also sample the restriction sample $\gamma$ with exponent $\alpha$, via a Poisson point process of Brownian excursions attached to $\R_-$, as explained in \cite {WW2005}.

Suppose now that $A\in\LA$, and define $H=H_A$ to be the unbounded connected component of $\HH \setminus A$ as before. By definition of $\LA$, the negative half-line still belongs to $\partial H_A$. If we restrict the loop-soup and the Poisson point process of Brownian excursions to those that stay in $H_A$, the restriction properties of the corresponding intensity measures imply immediately that one gets a sample of the Brownian loop-soup with intensity $c$ in $H_A$, and a sample of the Poisson point process of Brownian excursions away from the negative half-line in $H_A$, with intensity $\alpha$. In particular, because of the conformal invariance of these two underlying measures, it follows that  these Poissonian samples have the same law as the image under $\Phi_A^{-1}$ of the original loop and excursion soups in $\HH$.

\begin{figure}[ht!]
\begin{center}
\includegraphics[width=0.7\textwidth]{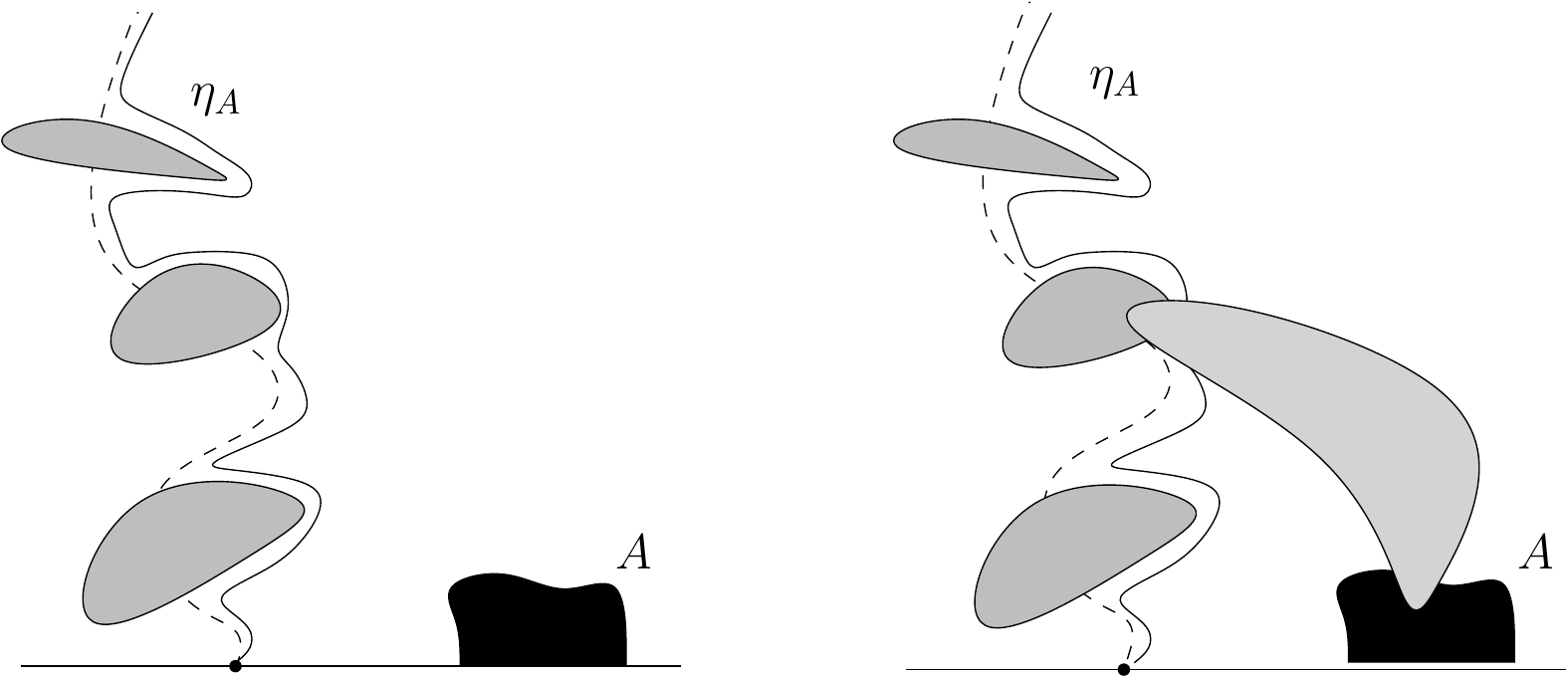}
\end{center}
\caption{\label{fig::radon_nikodym} $\eta=\eta_A$ if and only if there is no loop in $\Gamma$ that intersects $\eta_A$ and $A$.}
\end{figure}

Let us now first sample these items in $H_A$, and consider the right-most boundary of the curve $\eta_A$ defined just as $\eta$, but in $H_A$. Then, we sample those excursions and loops that do not stay in $H_A$, and we construct $\eta$ itself. One can note that either $\eta \not\subset H_A$ or $\eta=\eta_A$. Indeed, the only way in which $\eta$ can be different than $\eta_A$ is because of these additional loops/excursions, that do force $\eta$ to get out of $H_A$.
Hence, the event $\eta \subset H_A$ holds if and only if on the one hand the curve $\gamma$ stays in $H_A$ (recall that this happens with probability $\Phi_A'(0)^\alpha$), and on the other hand, no loop in the loop-soup does intersect both $\eta_A$ and $A$ (see Figure \ref{fig::radon_nikodym}).
It follows immediately that
for any $A\in\LA,$
\begin{equation*}
\frac{dP_{\HH}} {dP_{H_A}} (\eta) 1_{\eta \cap A = \emptyset}
=\Phi_A'(0)^{\alpha}\exp(-c L(\HH; A, \eta)) 1_{\eta \cap A = \emptyset}
\end{equation*}
where $L( \HH; A, \eta)$ denotes the mass (according to the Brownian loop-measure in $\HH$) of the set of loops that intersect both $A$ and $\eta$.

Equivalently,
\begin{equation}
\label {RNder}
\frac{dP_{H_A}} {dP_{\HH}} (\eta) 1_{\eta \cap A = \emptyset}
=1_{\eta \cap A = \emptyset} \Phi_A'(0)^{-\alpha}\exp(c L(\HH; A, \eta)).
\end{equation}
Note that this implies that
\begin {equation}
 \label {RNder2}
E_{\HH} \left( 1_{\eta \cap A = \emptyset} \exp( c L(\HH; A, \eta) \right) =
E_{H_A} (1_{\eta \cap A = \emptyset} \Phi_A'(0)^{\alpha})  = \Phi_A'(0)^{\alpha}
\end {equation}
(and the present argument in fact shows that the expectation in the left-hand side is actually finite).

We now wish to compare (\ref {RNder}) with features of SLE$_{\kappa}(\rho)$ processes. Let us now suppose that the curve $\eta$ is an SLE$_{\bar{\kappa}}(\bar{\rho})$ process
for some $\bar{\kappa}\le 4$ and $\bar{\rho} > -2$. We keep the same notations as in Subsection \ref {2.2}. For $A\in\LA,$ let $T$ be the (possibly infinite) first time
at which $\eta$ hits $A$.  For $t<T$, write $h_t := \Phi_{g_t(A)}$. Then (see \cite {JD2005}, Lemma 1), an It\^o formula calculation shows that
$$M_t=h_t'(W_t)^{a_1}h_t'(O_t)^{a_2}\left(\frac{h_t(W_t)-h_t(O_t)}{W_t-O_t}\right)^{a_3}\exp(\bar{c}L(\HH;A,\eta[0,t]))$$ is a local martingale (for $t < T$) where
$a_1=(6-\bar{\kappa})/(2\bar{\kappa})$, $a_2=\bar{\rho}(\bar{\rho}+4-\bar{\kappa})/(4\bar{\kappa})$, $a_3=\bar{\rho}/\bar{\kappa}$ and
and $\bar{c} = c(\bar{\kappa}) = (3\bar{\kappa}-8)(6-\bar{\kappa})/(2\bar{\kappa}) $ (note that such martingale calculations have been used on several occasions in related contexts, see e.g. \cite {JD2009} and the references therein).

It can be furthermore noted that $M_0=\Phi_A'(0)^{\bar{\alpha}}$ (and more generally, at those times when $O_t =W_t$, one puts $M_t = h_t'(W_t)^{\bar{\alpha}}  \exp(\bar{c}L(\HH;A,\eta[0,t]))$,
where $$\bar{\alpha} = \alpha(\bar{\kappa},\bar{\rho})= a_1+a_2+a_3=(\bar{\rho}+2)(\bar{\rho}+6-\bar{\kappa})/(4\bar{\kappa}).$$
One has to be a little bit careful, because (as opposed to the case where $\bar{\kappa} < 8/3$), $M_t$ is not bounded on $t<T$, so that we do not know if the local martingale stopped at $T$ is uniformly integrable (indeed the term involving $L( \HH; A , \eta [0,t])$ actually does blow up when $t \to T-$ and $T < \infty$).
However, even if some of the numbers $a_2$ and $a_3$ may be negative,
one always has (see \cite {JD2005}, the proof of Lemma 2-(i))
$$0 \le h_t'(W_t)^{a_1}h_t'(O_t)^{a_2}\left(\frac{h_t(W_t)-h_t(O_t)}{W_t-O_t}\right)^{a_3} \le 1.$$
Furthermore (see again \cite {JD2005}),  when $\eta \cap A = \emptyset$, then when $t \to \infty$, then $M_t$ converges to
$$ M_\infty : = \exp(\bar{c}L(\HH; A,\eta))$$
because each all the first three terms in the definition of $M_t$ converge to $1$. 

Note also that
$dM_t=M_tK_t \sqrt{\bar{\kappa}}dB_t$ where
$$K_t=a_1\frac{h_t''(W_t)}{h_t'(W_t)}+a_3\frac{h_t'(W_t)}{h_t(W_t)-h_t(O_t)}-a_3\frac{1}{W_t-O_t}.$$
Let $T_n$ denote the first (possibly infinite) time that the distance between the curve and $A$ reaches $1/n$. Then, for a fixed $A$, we see that $(M_{t \wedge T_n}, t \ge 0)$ is uniformly bounded by a finite constant. Hence, if $Q_{\HH}$ is the probability measure under which $W$ is the driving process of the SLE$_{\bar{\kappa}} (\bar{\rho})$ $\eta$ in $\HH$, we can define the probability measure $Q_n^*$ by $dQ_n^* / dQ_{\HH} = M_{T_n} / M_0$.
 Under $Q_n^*$, we have
$$dB_t=dB^*_t+K_tdt,\quad d h_t(W_t)=\sqrt{\bar{\kappa}}h'_t(W_t)dB_t^*+\frac{\bar{\rho}}{h_t(W_t)-h_t(O_t)}h'_t(W_t)^2dt.$$
This implies that $Q_n^*$ is the law of a (time-changed) SLE$_{\bar{\kappa}}(\bar{\rho})$ in $H_A$ up to the time $T_n$, which  happens to be the (possibly infinite) first time at which this curve gets to distance $1/n$ of $A$.

We can now note that by definition, the sequences $Q_n^*$ are compatible in $n$, so that there exists a probability measure $Q^*$ such that, under $Q^*,$ and for each $n$, the curve, up to time $T_n,$ is an SLE$_{\bar{\kappa}}(\bar{\rho})$ in $H_A$ up to the first time it is at distance $1/n$ of $A$.
But we also know that an SLE$_{\bar{\kappa}}(\bar{\rho})$ in $H_A$ almost surely does not hit $A$. Hence, $Q^*$ is just the law of SLE$_{\bar{\kappa}}(\bar{\rho})$ in $H_A$.

By the definition of $Q^*,$ we have that, for any  $n,$
$$\frac {dQ^*}{dQ_{\HH}} (\eta) 1_{d(\eta,A)\ge 1/n }=\frac {M_{T_n}}{M_0} 1_{d(\eta,A)\ge 1/n} = \frac {M_{\infty}}{M_0} 1_{d(\eta,A)\ge 1/n}.$$
Hence, we finally see that
$$\frac {dQ^*}{dQ_{\HH}} (\eta) 1_{d(\eta,A)> 0 }= \frac {M_\infty}{ M_0} 1_{d(\eta,A)> 0 } = \Phi_A'(0)^{- \bar{\alpha}} \exp ( \bar{c} L ( \HH; A , \eta)) 1_{ \eta \cap A = \emptyset}.$$
Comparing this with (\ref {RNder}), we conclude that $\bar{\kappa} = \kappa$ and that $\bar{\rho} = \tilde \rho (\kappa, \alpha)$.

\medbreak

Note that a by-product of this proof (keeping in mind that (\ref {RNder2}) holds) is that in fact the stopped martingale $M_{t \wedge T}$ is indeed uniformly integrable: It is a positive martingale such that
$$ E ( M_{T} ) = E ( \lim_{t \to \infty} M_{t \wedge T}) \ge E ( M_\infty 1_{ T = \infty }) = \Phi_A' (0)^\alpha =E( M_0).$$

\section {Proof of Lemma \ref {firstpart}}

We now describe the steps of the proof of Lemma \ref {firstpart}. Quite a number of these steps are almost identical to ideas developed in \cite {SW2010}.
We will therefore not always provide all details of those parts of the proof. Let us first note that the law of $\eta$ is obviously scale-invariant, and that we already have seen that it is almost surely a simple curve.
Furthermore, we know (for instance using the construction of $\gamma$ via a Poisson point process of Brownian excursions, or via its SLE$_{8/3}(\rho)$ description), that almost surely, the Lebesgue measure of $\gamma \cap (-\infty, 0)$ is zero. By construction (since $\eta \cap (-\infty, 0)$ is a subset of this set), the Lebesgue measure of $\eta \cap (-\infty, 0)$ is also $0$. Hence, in order to prove the lemma, it only remains to check the ``conformal Markov'' property i.e. the last item in  the characterization of  SLE$_{\kappa}(\rho)$ processes derived in Subsection \ref {2.2}.

\subsection {Straight exploration and the pinned path}

A first idea will be not to focus only on the curve $\eta$, but to also keep track of the CLE loops that lie to its right. In other words, we will consider
half-plane configurations $(\eta,\Lambda)$, where -- as before --  $\eta$ is a curve in $\overline{\HH}$ from 0 to $\infty$ that does not touch $(0,\infty)$ and $\Lambda$ is a
loop configuration in the connected component of $\HH\setminus \eta$ that has $(0,\infty)$ on its boundary (we say that it is the connected component to the right of $\eta$).
The conformal restriction property of the CLE shows that the following two constructions are equivalent:

\begin {itemize}
 \item Construct $\eta$ as in the statement (via a CLE $\Gamma$ and a restriction path $\gamma$), and consider $\Lambda$ to be the collection of loops in the CLE $\Gamma$ (that one used to construct $\eta$) that lie to the right of $\eta$.
 \item First sample $\eta$, and then in the connected component $H_\eta$ of $\HH \setminus \eta$ that lies to the right of $\eta$, sample an independent CLE that we call $\Lambda$.
\end {itemize}

It turns out that the couple $(\eta, \Lambda)$ does satisfy a simple ``restriction-type'' property, that one can sum up as follows: For a given $A \in \LA$, let us condition on the event $\{ \eta \cap A = \emptyset \}$. Then, one can define the collection $\tilde \Lambda_A$ of loops of $\Lambda$ that intersect $A$, and the unbounded connected component $\tilde H_A$  of $\HH \setminus (A \cup \tilde \Lambda_A)$. We also denote by $\Lambda_A$ to be the collection of loops of $\Lambda$ that stay in $\tilde H_A$.  Let $\Psi= \Psi (\tilde \Lambda_A, A)$ denote the conformal map from $\tilde H_A$ onto $\HH$ with $\Psi (0)= 0$ and $\Psi (z) \sim z$ when $z \to \infty$. Then, the conditional law of
$( \Psi ( \eta) , \Psi ( \Lambda_A ) )$ (conditionally on $\eta\cap A=\emptyset$) is identical to the original law of $(\eta, \Lambda)$. This is a direct consequence of the construction of $(\eta, \Lambda)$ and the restriction properties of $\gamma$ and $\Gamma$.

\begin{figure}[ht!]
\begin{center}
\includegraphics[width=0.8\textwidth]{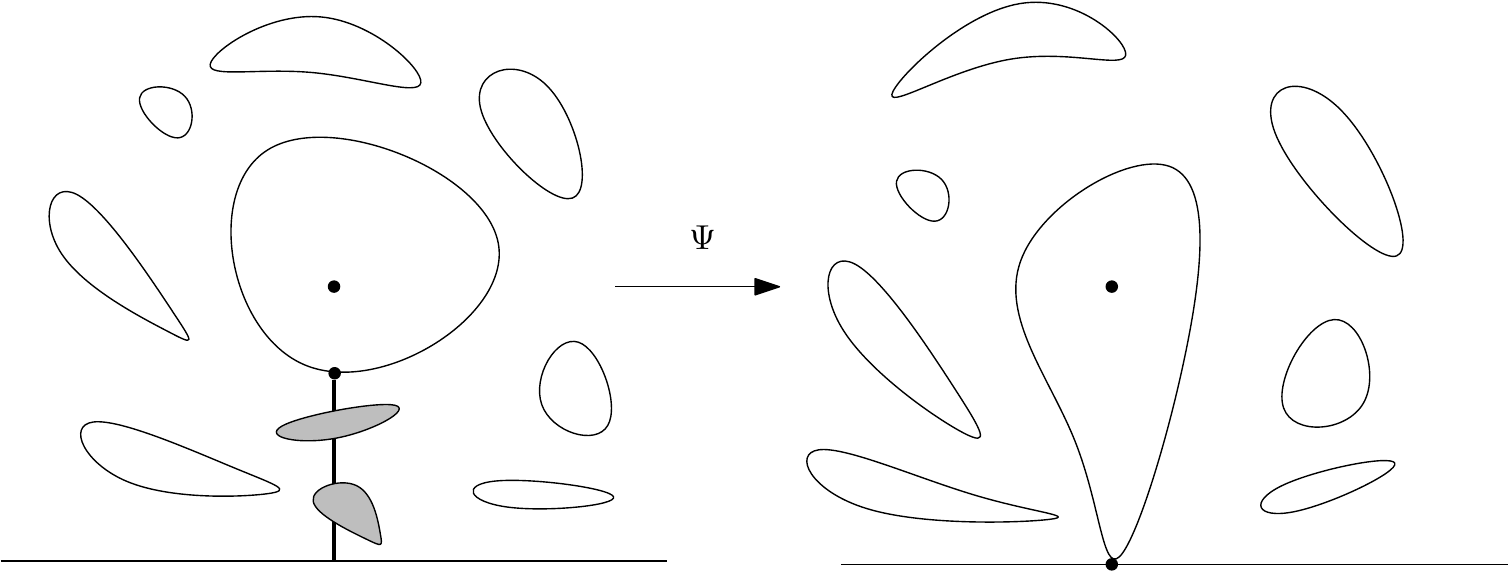}
\end{center}
\caption{\label{fig::exploration_CLE} Discovering the loop that surrounds $i$ in a CLE defines a pinned loop (see \cite {SW2010})}
\end{figure}

This restriction property is of course reminiscent of the restriction property of CLEs themselves. In \cite {SW2010}, the restriction property of CLE was exploited
as follows: Fix one point in $\HH$ (say the point $i$) and discover all loops of the CLE that lie on the segment $[0,i]$ (by moving upwards on this segment) until one discovers the loop that surrounds $i$ (see Figure \ref{fig::exploration_CLE}).
This can be approximated by iterating discrete small cuts, discovering the loops that interesect these cuts and repeating the procedure. The outcome was a description of the law of the loop that surrounds $i$ at the ``moment'' at which one discovers it (see Proposition 4.1 in  \cite {SW2010}).

Here, we use the very same idea, except that the goal is to cut in the domain until one reaches the curve $\eta$
(note that in the CLE case, the marked point $i$ is an interior point of $\HH$ and that here, the marked points $0$ and $\infty$ on the boundary do also correspond to the choice of two degrees of freedom in the conformal map).
 We can for instance do this by moving upwards on the vertical half-line $L:= 1 + i \R_+$;
a simple $0$-$1$ law argument shows that almost surely, the curve $\gamma$ does intersect $L$, and that therefore $\eta \cap L \not= \emptyset$ too. Let $\eta_T$ denote the point of $\eta \cap L$ with smallest $y$-coordinate. One way to find it, is to move on $L$ upwards until one meets $\eta$ for the first time. This can be approximated also by ``exploration steps'', in a way that is almost identical to the explorations of CLEs described in \cite {SW2010}. We refer to that paper for rather lengthy details, the arguments really just mimic those to that paper. The conclusion, analogous to Proposition 4.1 in \cite {SW2010} is that (see Figure \ref{fig::exploration_SLE}):

\begin{figure}[ht!]
\begin{center}
\includegraphics[width=0.9\textwidth]{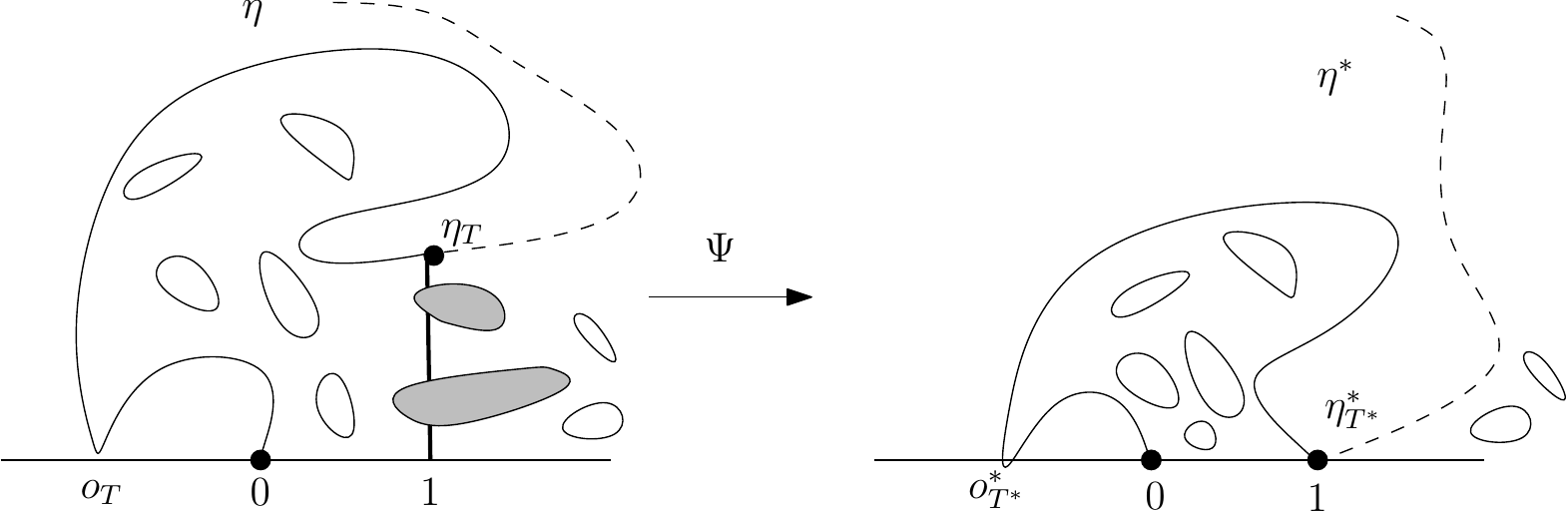}
\end{center}
\caption{\label{fig::exploration_SLE}Discovering $\eta$ in half-plane configuration defines a pinned path}
\end{figure}

\begin {lemma}
\label {straightexplo}
 The conditional law of $\eta$ conditionally on the event that $\eta$ passes through the $\epsilon$-neighborhood of $1$, converges as $\epsilon \to 0$ to the distribution of
 $\eta^*:= \Psi (\eta)$, where $\Psi$ is the conformal map from $\tilde H_{[1, \eta_T]}$ onto $\HH$ that maps the triplet $( 0, \eta_T, \infty)$ onto $(0, 1, \infty)$.
\end {lemma}
We will call  $\eta^*$  a ``pinned'' path, as in \cite {SW2010}. Note that this construction also shows that $\eta^*$ is independent of $\Psi$.

\subsection {Restriction property for the pinned path}

When $\eta^*$ is such a pinned path, then $\HH \setminus \eta^*$ has several connected components, and we call $U_0$ the connected components with $(0,1)$ on its boundary and $U_+$  the one with $(1,\infty)$ on its boundary (see Figure \ref{fig::indep_future_past_2}). If one first samples $\eta^*$ and then in $U_0$ and $U_+$ samples two independent CLE$_\kappa$'s , then one gets a ``pinned configuration'' $(\eta^*, \Lambda^*)$.

This pinned configuration inherits the following restriction property from $(\eta, \Lambda)$: Suppose that $A \in \LA$ with $d(1, A) > 0$, and condition on $A \cap \eta^* = \emptyset$. Then, define
$H_A^*$ for $(\eta^*,\Lambda^*)$ just as $\tilde H_A$ in the case of $(\eta, \Lambda)$. Note that $0$ and $1$ are both boundary points of $H_A^*$ so that it is possible to define the conformal transformation
$\Phi_A^*$ from $H_A^*$ onto $\HH$ that fixes the three boundary points $0$, $1$ and $\infty$.

Then, the conditional law of $\Phi_A^* (\eta^*)$ (conditionally on $\eta^* \cap A = \emptyset$) is equal to the law of $\eta^*$ itself. This result just follows by passing to the limit the restriction property of $(\eta, \Lambda)$.

\begin{figure}[ht!]
\begin{center}
\includegraphics[width=0.9\textwidth]{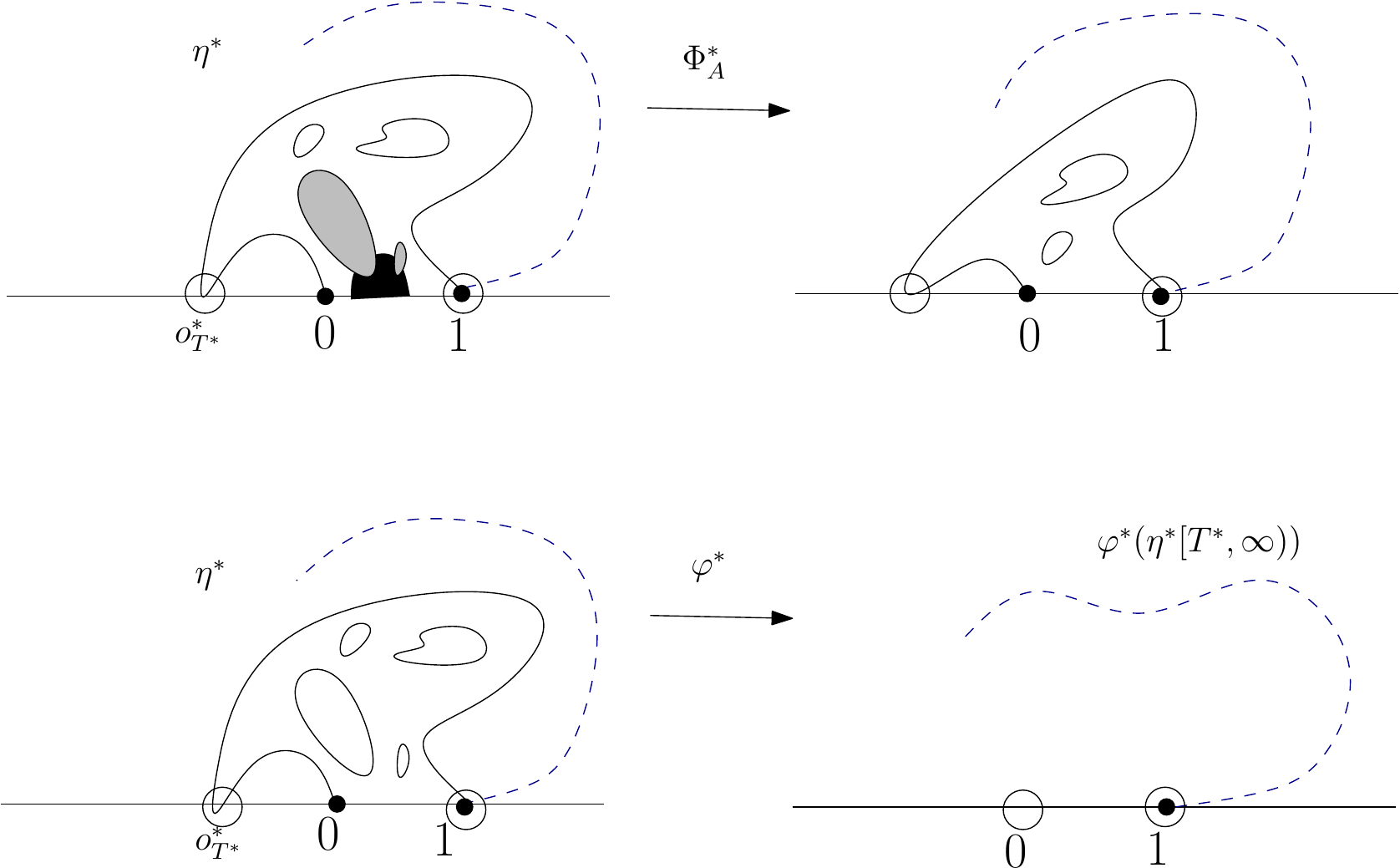}
\end{center}
\caption{\label{fig::indep_future_past_1} Definitions of $\Phi_A^*$ and  $\varphi^*$}
\end{figure}

Let us define $T^*$ the time at which $\eta^*_{T^*}= 1$, and $o^*_{T^*}$ as the leftmost point in  $\eta^*[0,T^*]\cap\R_-$ (note that depending on the value of $\rho$, it may be the case that $o^*_{T^*}=0$). Denote by $\varphi^*$ the conformal map from the unbounded connected component of
$\HH \setminus {\eta^* [0, T^*]}$ onto $\HH$, that maps the triplet $(o^*_{T^*},1,\infty)$ onto $(0,1,\infty)$ (see Figure \ref{fig::indep_future_past_1}).

Let us now consider a set $A \in \LA$ that is also at positive distance from $[1, \infty)$, i.e. that is attached to the segment $[0,1]$ (we call $\LA_{[0,1]}$ this set of events).
Then, the following restriction property will be inherited from the restriction property of $(\eta^*, \Lambda^*)$:

\begin {lemma}
\label {indep}
 The curve $\varphi^* ( \eta^* [ T^* , \infty ))$ is independent of the event $\eta^* [0, T^*] \cap A = \emptyset$.
\end {lemma}

Indeed, if one conditions on the event $\eta^* [0, T^*] \cap A = \emptyset$ (which is the same as $\eta^*\cap A=\emptyset$), then the conditional law of
$\Phi_A^* ( \eta^*)$ is that of $\eta^*$ itself, so that $\eta^* [0, T^*] \cap A = \emptyset$ and $\Phi_A^* ( \eta^*)$ are independent.
But $\varphi^* ( \eta^* [ T^* , \infty ))$ can also be recovered from $\Phi_A^* ( \eta^*)$ (see Figure \ref{fig::indep_future_past_1}). This implies the Lemma.

A direct consequence of the lemma is therefore that $\eta^* [0, T^*]$ and $\varphi^* ( \eta^* [ T^* , \infty ))$ are independent. Indeed, the $\sigma$-field generated by the family of events of the type
 $\eta^* [0, T^*] \cap A = \emptyset$ when $A \in \LA_{[0,1]}$ (which is stable by finite intersections) is exactly $\sigma ( \eta^* [0, T^*])$.

\subsection {General explorations and consequences}

\begin{figure}[ht!]
\begin{center}
\includegraphics[width=0.9\textwidth]{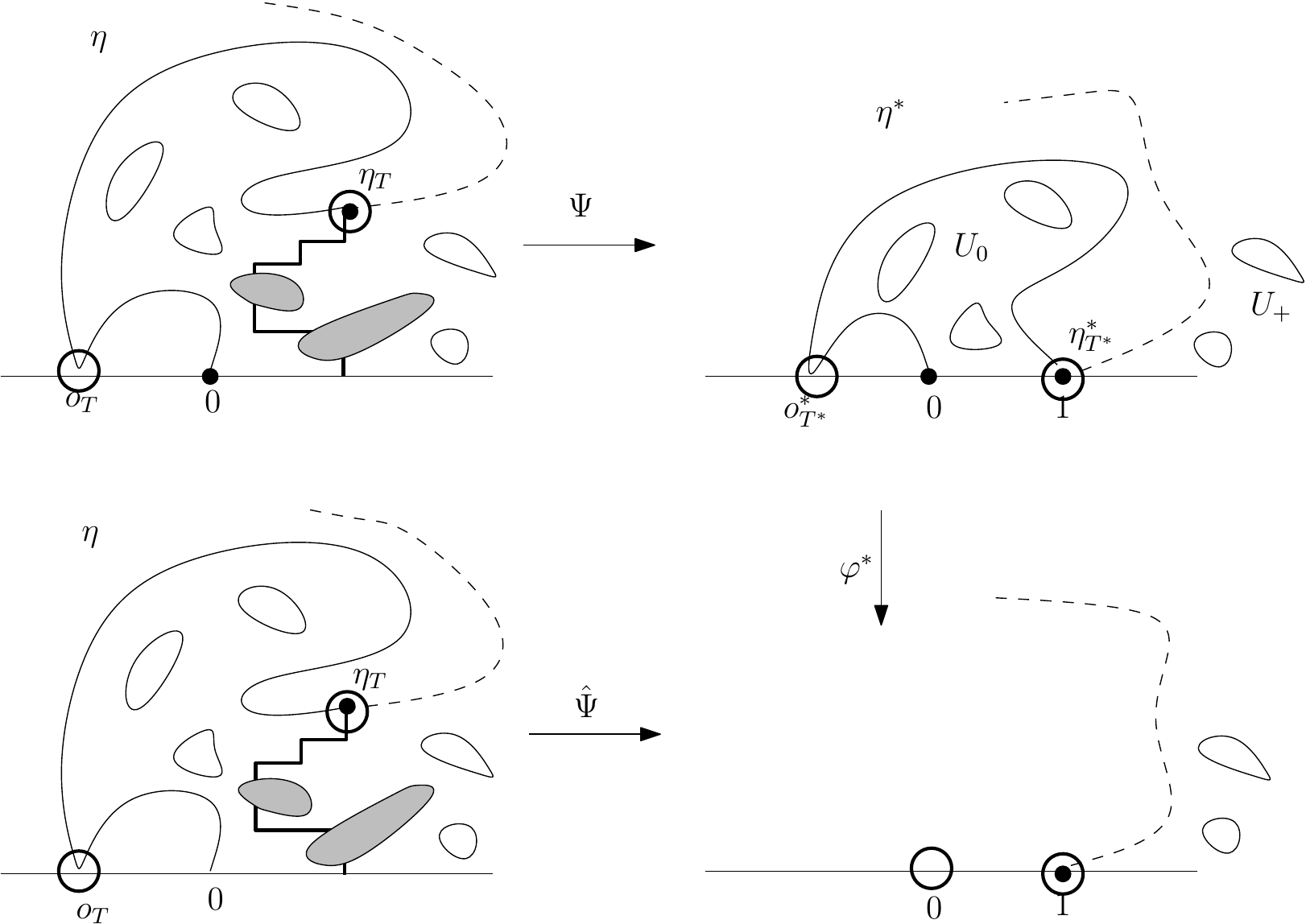}
\end{center}
\caption{\label{fig::indep_future_past_2} $\Psi$, $\varphi^*$  and $\hat{\Psi}=\varphi^*\circ\Psi$.}
\end{figure}

In fact, just as in \cite {SW2010}, it is easy to see that the argument that leads to Lemma \ref {straightexplo} can be generalized to other curves than the straight line $L$.
In particular, if we choose $L$ to be any oriented simple curve on the grid $\delta (\Z \times \N)$ that starts on the positive half-line and disconnects $0$ from infinity in $\HH$, then define $\eta_T$ to be the point of $\eta$ that $L$ meets first, and let $\tilde L$ denote the part of $L$ until it hits $\eta_T$. We then define $\tilde H$ as the unbounded connected component of the set obtained by removing from $\HH \setminus \tilde L$ all the loops of $\Lambda$ that intersect $\tilde L$. Let $\Psi$ denote the conformal map from $\tilde H$ onto $\HH$ that sends the triplet $(0,\eta_T,\infty)$ onto $(0,1,\infty).$ Let $\hat{H}$ be the unbounded connected component of the set obtained by removing from $\HH$ the union of $\eta[0,T],\tilde L$ and the loops in $\Lambda$ that intersect $\tilde L.$ Let $\hat{\Psi}_{\tilde{L}}$ denote the conformal map from $\hat{H}$ onto $\HH$ that sends the triplet $(o_T,\eta_T,\infty)$ onto $(0,1,\infty)$ (see Figure \ref{fig::indep_future_past_2}).
Then the same arguments than
the ones used to derive
Lemma \ref{straightexplo} imply that $\Psi(\eta)$ has the same law as pinned path $\eta^*.$ Combined with Lemma \ref{indep}, this implies that $\hat{\Psi}_{\tilde{L}}(\eta[T,\infty))$ is independent of $\eta[0,T].$

\begin{figure}[ht!]
\begin{center}
\includegraphics[width=0.9\textwidth]{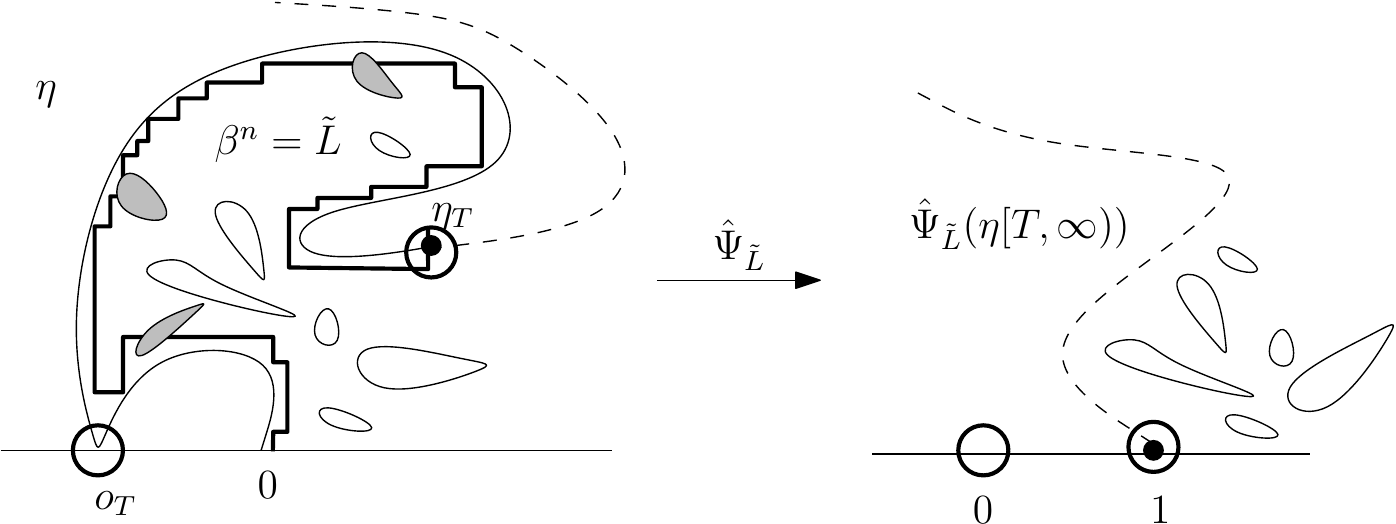}
\end{center}
\caption{\label{fig::approximation} $\hat{\Psi}_{\tilde{L}}$ maps the triplet $(o_T,\eta_T,\infty)$ onto $(0,1,\infty).$}
\end{figure}

The next step of the proof is again almost identical to the corresponding one in \cite {SW2010}:
Fix a time $T$ and suppose that $\eta_T\not\in\R.$ Let $\beta^n$ be an approximation of $\eta[0,T]$ from right on the lattice $2^{-n}(\Z\times\N)$ (see Figure \ref{fig::approximation}). Then for any deterministic piecewise linear path $\tilde{L},$ on the event $\{\beta^n=\tilde{L}\}$, the probability that $\tilde{L}$ intersects some macroscopic loop in $\Lambda$ is very small when $n$ is large enough, so that $\hat{\Psi}_{\tilde{L}}(\eta[T,\infty))$ is very close to $f_T(\eta[T,\infty))$ on this event. Since $\hat{\Psi}_{\tilde{L}}(\eta[T,\infty))$ is independent of $\eta[0,T],$ by passing to the limit (as $n\to\infty$), we get that $f_T(\eta[T,\infty))$ is independent of $\eta[0,T]$ as desired. This is exactly the conformal Markov property that was needed to conclude the proof of Lemma \ref {firstpart}.

\section {Consequences for second-moment estimates}

In order to illustrate how the present construction can be used in order to derive directly some properties of SLE$_{\kappa}(\rho)$ processes, we are going to derive in this section
some information about the intersection of SLE$_{\kappa}(\rho)$ processes and the real line. Analogous ideas have been used in \cite {NW2011} to study the dimension of the CLE gasket, but the situation here is even more convenient.

Recall that from the definition, we know that the SLE$_{\kappa}(\rho)$ process $\eta$, from $0$ to $\infty$ in $\HH$ does not touch the positive half-line, but -- as we already mentioned --, its definition via the Loewner equation and  Bessel processes shows that it touches almost surely the negative half-line as soon as $ \rho < (\kappa/2) -2 $. For instance, for $\kappa =4$, this will happen for $\rho \in (-2, 0)$, while for $\kappa = 3$, this will occur for $\rho \in (-2, -1/2)$. Here for obvious reasons, we will restrict ourselves to the case where $\kappa \in (8/3, 4]$.

\begin{proposition}\label{dimension}
For $\kappa\in(8/3,4]$ and $\rho\in (-2, -2 + \kappa/2)$, then the Hausdorff dimension of $\eta \cap \R_-$ is almost surely equal to
$1-(\rho+2)(\rho+4-\kappa/2)/\kappa$.
\end{proposition}
Note that this result is also derived in \cite{MW2012} for all $\kappa\in(0,8)$ and $\rho \in (-2, -2 + (\kappa/2))$ using the properties of flow lines of GFF introduced in \cite{MS2012a}.

\medbreak

Before turning our attention to the proof of this result, let us first focus on the following related question: Let us fix $c \in (0,1)$ and $\alpha > 0$.
Consider on the one hand a Brownian loop-soup with intensity $c$ in the upper half-plane, and its corresponding CLE$_\kappa$ sample consisting of the outermost boundaries of the loop-soup clusters, as in \cite {SW2010}.

On the other hand, consider a Poisson point process $(b_j, j \in J)$ of Brownian excursions away from the real line in $\HH$, with intensity $\alpha$.
Each of these excursions $b_j$ has a starting point $S_j$ and an endpoint $E_j$ that both lie on the real axis.

\begin{figure}[ht!]
\begin{center}
\includegraphics[width=0.7\textwidth]{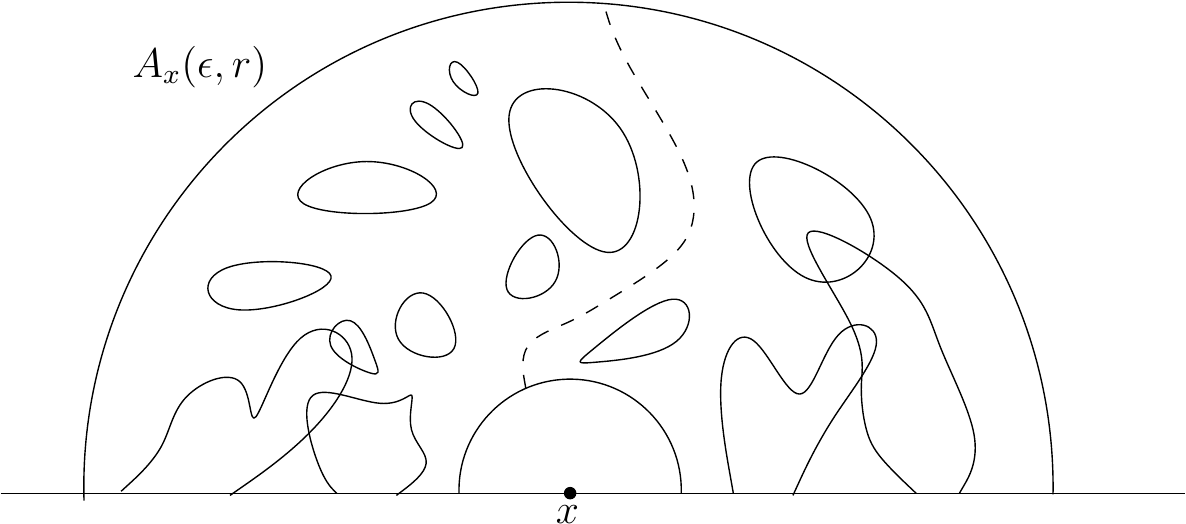}
\end{center}
\caption{\label{fig::annulus} Event $E_x(\eps,r)$: $x$ is not disconnected from $\infty$ by the excursions and loops.}
\end{figure}

For each point $x$ on the real line, for each $\epsilon < r $, we define the semi-ring
$$A_x (\epsilon , r ) : = \{ z \in \overline \HH \ : \ \epsilon < | z - x | < r \}.$$
For each given $\epsilon$ and $r$, we can artificially restrict ourselves to those Brownian loops and excursions that stay in $A_x (\epsilon, r)$. We define the event
$E_x ( \epsilon, r )$ that the union of all these paths does not disconnect $x$ from infinity in $\HH$ (see Figure \ref{fig::annulus}).

Clearly, the probability of $E_x (\epsilon, r)$ is in fact a function of $\epsilon/r$ and does not depend on $x$. Let us denote this probability by $p(\epsilon/r)$.
It is elementary to see that for all $\epsilon, \epsilon' < 1$, $$p(\epsilon \epsilon') \le p(\epsilon) p(\epsilon').$$ Indeed, if one divides $A_0 (\epsilon \epsilon', 1)$ into the two semi-annuli $A_0 (\epsilon \epsilon' , \epsilon)$ and $A_0 (\epsilon, 1)$, one notices that
$$ E_0 (\epsilon \epsilon , 1 ) \subset E_0 (\epsilon \epsilon', \epsilon) \cap E_0 ( \epsilon, 1)$$
and the latter two events are independent, due to their Poissonian definition.

On the other hand, for some universal constant $C$, we know that for all $\epsilon, \epsilon' < 1/4$,
\begin {equation}
 \label {superadd}
p(8 \epsilon \epsilon') \ge C p(\epsilon) p(\epsilon').
\end {equation}
Indeed, let us consider the following three events:
\begin{itemize}
\item
$U_1$: No CLE loop touches both $\{ z \ : \ |z| = 2 \}$ and $\{ z \ : \ |z|= 4 \}$
\item
$U_2$: No Brownian excursion touches both $\{ z \ : \ |z|=1 \}$ and $\{ z \ : \  |z| = 2 \}$.
\item
$U_3$: No Brownian excursion touches both $\{ z \ : \ |z|=4 \}$ and $\{ z \ : \ |z|=8 \}$.
\end {itemize}
All the events $U_1$, $U_2$, $U_3$, $E_0 (8\epsilon, 8)$ and $E_0 (1, 1/\epsilon')$ are decreasing events of the Poisson point processes of loops and excursions
(i.e. if an event fails to be true, then adding an extra excursion or loop will not fix it). Hence, they are positively correlated.
Furthermore, we have chosen these events in such a way that
$$ \left( U_1 \cap U_2 \cap U_3 \cap E_0 (8\epsilon, 8) \cap E_0 (1, 1/\epsilon' )\right)  \subset E_0 ( 8 \epsilon, 1/ \epsilon').$$
The fact that $c \le 1$ ensures that the events $U_1$, $U_2$ and $U_3$ have a positive probability. Putting the pieces together, we get that
$$ p( 8 \epsilon \epsilon' ) = P (  E_0 ( 8 \epsilon, 1/ \epsilon')) \ge P (U_1 \cap U_2 \cap U_3) p(\epsilon) p(\epsilon')$$
from which (\ref {superadd}) follows.
Hence, if we define $q(\epsilon):= C p( 8 \epsilon)$, we get $q ( \epsilon \epsilon') \ge q ( \epsilon) q ( \epsilon')$.

These properties of $p(\epsilon)$ and $q (\epsilon)$ ensure that there exists a positive finite $\beta$ and a constant $C'$ such that for all $\epsilon < 1/8$,
$$ \epsilon^{\beta} \le p( \epsilon ) \le C' \epsilon^\beta.$$

\medbreak

Let us now focus on the proof of the proposition. First, let us note that a simple $0$-$1$ argument (because the studied property is invariant under scaling) shows that
there exists $D$ such that almost surely, the dimension of $\eta \cap \R_-$ is equal to $D$. Furthermore, we can use scale-invariance again to see that in order to prove that $D$ is equal to some given value $d$,
it suffices to prove that on the one hand, almost surely, the Hausdorff dimension of $\eta \cap [-2, -1]$ does not exceed $d$,  and that on the other hand, with positive probability, the Hausdorff dimension of $\eta \cap [-2, -1]$ is equal to $d$.

Let us now note that if a point $x \in [-2, -1]$ belongs to the $\epsilon$-neighborhood  $K_\epsilon$ of
$\eta$, then it implies that $E_x (\epsilon, 1)$ holds.
Hence, the first moment estimate implies readily that almost surely, the Minkovski dimension of $\eta \cap [-2, -1]$ is not greater than  $ 1 - \beta$, and therefore that $D \le 1 - \beta$.

\medbreak
In order to prove that with positive probability, the dimension of $\eta \cap [-2, -1]$ is actually equal to $1-  \beta$, we can make the following two observations.
\begin {itemize}
 \item Suppose that $x \in [-2, -1]$ and that $E_x ( \epsilon/2, 8)$ holds. Suppose furthermore that no excursion in the Poisson point process of excursions attached to $(- \infty, -6)$ does intersect the ball of radius $4$ around the origin, no excursion in the Poisson point process excursions attached to $(-2,0)$ exits the ball of radius $4$ around $-2$. Suppose furthermore that no loop in the CLE (in $\HH$) intersects both the circle of radius 4 and 6 around the origin.
 Note that these two events have positive probability and are positively correlated to $E_x (\epsilon /2 , 8)$ (they are all decreasing events of the Poisson point processes of loops and excursions). Then, by construction, $x$ is necessarily in the $\epsilon$-neighborhood of $\eta$. It therefore follows that for some constant $c'$, for all $x \in [-2, -1]$,
 $$ P ( x \in K_\epsilon ) \ge c' \epsilon^{\beta}.$$

 \item Suppose now that $-2< x < y < -1$, that $y-x < 1/4$ and that $\epsilon < (y-x)/4$. Clearly, if both $x$ and $y$ belong to $K_\epsilon$, then it means that the three events
 $E_x ( \epsilon, (y-x)/2)$, $E_y (\epsilon, (y-x)/2)$ and $E_x (2(y-x) , 1/2)$ hold. These three events are independent, and the previous estimates therefore yield that there exists a constant $c''$ such that
 $$ P ( x \in K_\epsilon, y \in K_\epsilon) \le c'' \frac {\epsilon^{2 \beta}}{ (y-x)^\beta}.$$
\end {itemize}
Standard arguments (see for instance \cite {MP}) then imply that with positive probability, the dimension of $\eta \cap [-2, -1]$ is not smaller than $1 - \beta$.
This concludes the proof of the fact that almost surely, the Hausdorff dimension of $\eta \cap (-\infty, 0) $ is almost surely equal to $1-\beta$.

\medbreak
In order to conclude, it remains to compute the actual value of $\beta$. A proof of this is provided in \cite {MW2012} using the framework of flow lines of the Gaussian Free Field. Let us give here an outline of how to compute $\beta$ bypassing the use of the Gaussian Free Field, using the more classical direct way to derive the values of such exponents i.e. to exhibit a fairly simple martingales involving the derivatives of the conformal maps at a point, and then to use this to estimate the probability that the path ever reaches a small distance of this point:
Consider the SLE$_{\kappa}(\rho)$ process in $\HH$ from 0 to $\infty,$ and keep the same notations as in Subsection \ref {2.2}. First, one can note that for any real $v$,
$$M_t=g_t'(-1)^{v(\kappa v +4-\kappa)/4}(W_t-g_t(-1))^{v}(O_t-g_t(-1))^{v\rho/2}$$
is a local martingale. We then choose $v=(\kappa-8-2\rho)/\kappa,$ and define $\tilde \beta:=(\rho+2)(\rho+4-\kappa/2)/\kappa$ as well as
$$\Upsilon_t=\frac{O_t-g_t(-1)}{g_t'(-1)},\quad N_t=\frac{O_t-g_t(-1)}{W_t-g_t(-1)},\quad \tau_{\eps}=\inf\{t: \Upsilon_t=\eps\}.$$
Then $M_t=\Upsilon_t^{-\tilde \beta}N_t^{-v}$. Furthermore, the probability that the curve gets within the ball centered at $-1$ of radius $\eps$ is
comparable to $P(\tau_{\eps}<\infty).$ But, one has
$$P(\tau_{\eps}<\infty)=E(M_{\tau_{\eps}}N_{\tau_{\eps}}^v 1_{\tau_{\eps}<\infty})\eps^{\tilde \beta}=E^*(N_{\tau_{\eps}}^v)\eps^{\tilde \beta}$$
where $P^*$ is the measure $P$ weighted by the martingale $M.$ Under $P^*,$ we have that $\tau_{\eps}<\infty$ almost surely and that $E^*(N_{\tau_{\eps}}^v)$ is bounded both sides by universal constants independent of $\eps.$ It follows that indeed $\beta=\tilde \beta$.

\medbreak

We conclude with the following two remarks: 
\begin {itemize}
 \item Similar second-moment estimates can be performed for other questions related to SLE$_\kappa (\rho)$ processes for $\kappa \in (8/3, 4]$ and $\rho > -2$. 
 For instance the boundary proximity estimates from Schramm and Zhou \cite {SZ} can be generalized/adapted to the SLE$_\kappa (\rho)$ cases. We leave this to the interested reader.
 \item 
 It is proved in \cite {MS2012a} that the left boundary of an SLE$_{\kappa_0} (\rho_0)$ process for $\kappa_0 > 4$ and $\rho_0 > -2$ is an SLE$_{\kappa_1} ( \rho_1, \rho_2)$ process 
 for $\kappa_1 = 16 / \kappa_0$ with an explicit expression of $\rho_1$ and $\rho_2$ in terms of $( \kappa_0, \rho_0)$ (this is the ``generalized SLE duality'').
 Hence, it follows from Proposition \ref {dimension} that its statement (i.e. the formula for the Hausdorff dimension) in fact holds true for all $\kappa \in (4,6)$ as well. However, since the Gaussian Free Field approach is anyway used in the derivation of this generalized duality result, it is 
 rather natural to use also the Gaussian Free Field in order to derive the second moments estimates, as done in \cite {MW2012}.
 The same remark applies to the intersection of the right boundary of an SLE$_{\kappa_0} (\rho_0)$ when $\kappa_0 > 4$ and $\rho_0 \in (-2, 0)$; the Hausdorff 
 dimension of the intersection of this right boundary with $\R_-$ then turns out to be 
 $$1 - \frac {(\rho_0+2)(\rho_0 +  (\kappa_0/2))}{ \kappa_0} = - \rho_0 \left( \frac {\rho_0 + 2}{\kappa_0}  + \frac 1 2 \right).$$

\end {itemize}

\bigbreak

{D\'epartement de Math\'ematiques, Universit\'e Paris-Sud, 91405 Orsay Cedex, France\\

wendelin.werner@math.u-psud.fr \\

hao.wu@math.u-psud.fr\\


\begin{thebibliography}{99}

\bibitem{CDHKS}
Dmitry Chelkak, Hugo Duminil-Copin, Cl\'ement Hongler, Antti Kemppainen and Stanislav Smirnov.
\newblock Convergence of Ising interfaces to Schramm’s SLEs.
\newblock preliminary version, 2012

\bibitem{CS2012}
Dmitry Chelkak and Stanislav Smirnov.
\newblock Universality in the 2{D} {I}sing model and conformal invariance of
  fermionic observables.
\newblock {\em Invent. Math.}, 189:515--580, 2012.

\bibitem{JD2005}
Julien Dub{\'e}dat.
\newblock {${\rm SLE}(\kappa,\rho)$} martingales and duality.
\newblock {\em Ann. Probab.}, 33:223--243, 2005.

\bibitem{JD2009}
Julien Dub\'edat.
\newblock {Duality of Schramm-Loewner Evolutions.}
\newblock {Ann. Sci. ENS.}, 42:697--724, 2009.

\bibitem{LSW2003}
Gregory Lawler, Oded Schramm, and Wendelin Werner.
\newblock Conformal restriction: the chordal case.
\newblock {\em J. Amer. Math. Soc.}, 16:917--955, 2003.

\bibitem{GL2005}
Gregory~F. Lawler.
\newblock {\em Conformally invariant processes in the plane}, volume 114 of
  {\em Mathematical Surveys and Monographs}.
\newblock American Mathematical Society, Providence, RI, 2005.

\bibitem{LW2004}
Gregory~F. Lawler and Wendelin Werner.
\newblock The {B}rownian loop soup.
\newblock {\em Probab. Theory Related Fields}, 128:565--588, 2004.

\bibitem{MP}
Peter M\"orters and Yuval Peres.
\newblock Brownian motion.
\newblock Cambridge University Press, 2010.

\bibitem{MS2012a}
Jason Miller and Scott Sheffield.
\newblock Imaginary geometry i: Interacting SLEs.
\newblock Preprint, 2012.

\bibitem{MS2012b}
Jason Miller and Scott Sheffield.
\newblock Imaginary geometry ii: reversibility of SLE$_{\kappa}(\rho_1;\rho_2)$
  for $\kappa\in (0,4)$.
  \newblock Preprint, 2012.

\bibitem{MW2012}
Jason Miller and Hao Wu.
\newblock In preparation, 2012.

\bibitem{NW2011}
{\c{S}}erban Nacu and Wendelin Werner.
\newblock Random soups, carpets and fractal dimensions.
\newblock {\em J. Lond. Math. Soc. (2)}, 83:789--809, 2011.

\bibitem{RY1999}
Daniel Revuz and Marc Yor.
\newblock {\em Continuous martingales and {B}rownian motion}.
\newblock Springer-Verlag, Berlin, third edition, 1999.

\bibitem{OS2000}
Oded Schramm.
\newblock Scaling limits of loop-erased random walks and uniform spanning
  trees.
\newblock {\em Israel J. Math.}, 118:221--288, 2000.

\bibitem{SSW}
Oded Schramm, Scott Sheffield and David Wilson.
\newblock {Conformal radii
for conformal loop ensembles.}
\newblock {\em Comm. Math. Phys.}, 288:43--53, 2009.


\bibitem{SchrammW2005}
Oded Schramm and David~B. Wilson.
\newblock S{LE} coordinate changes.
\newblock {\em New York J. Math.}, 11:659--669, 2005.

\bibitem{SZ}
{Oded Schramm and Wang Zhou.}
\newblock Boundary proximity of SLE.
\newblock {\em Probab. Theor.  Rel. Fields.} 146:435--450, 2010.

\bibitem {Sheffield}
Scott Sheffield.
\newblock Exploration trees and conformal loop ensembles.
\newblock {\em Duke Math. J.},   147:79-129, 2009.


\bibitem{SW2010}
Scott Sheffield and Wendelin Werner.
\newblock Conformal loop ensembles: The markovian characterization and the
  loop-soup construction.
  \newblock  {\em Ann. Math.}, 176:1827-1917, 2012.

\bibitem{WW2003}
Wendelin Werner.
\newblock S{LE}s as boundaries of clusters of {B}rownian loops.
\newblock {\em C. R. Math. Acad. Sci. Paris}, 337:481--486, 2003.

\bibitem{WW2004}
Wendelin Werner.
\newblock Girsanov's transformation for {${\rm SLE}(\kappa,\rho)$} processes,
  intersection exponents and hiding exponents.
\newblock {\em Ann. Fac. Sci. Toulouse Math. (6)}, 13:121--147, 2004.

\bibitem{WW2005}
Wendelin Werner.
\newblock Conformal restriction and related questions.
\newblock {\em Probab. Surv.}, 2:145--190, 2005.

\bibitem {WW}
Wendelin Werner and Hao Wu.
\newblock On conformally invariant CLE explorations.
\newblock preprint, 2011.

\bibitem {zhan0}
Dapeng Zhan.
\newblock Reversibility of chordal SLE.
\newblock {\em Ann. Probab.}, 36:1472--1494, 2008.
\bibitem{zhan}
Dapeng Zhan.
\newblock Reversibility of some chordal SLE($\kappa;\rho$) traces.
\newblock {\em J. Stat. Phys.}, 139:1013-1032, 2010.

\end{thebibliography}
\end{document}